\documentclass[letterpaper, 10 pt, conference]{IEEEconf} 

\IEEEoverridecommandlockouts 

\usepackage{cite}
\usepackage{amsmath}
\usepackage{amsfonts}
\usepackage{color}
\usepackage{graphicx}
\graphicspath{ {./images/} }

\usepackage{algorithm}
\usepackage{algpseudocode}

\newcommand{\alias}{C:/Users/yuzhe/Dropbox/bib/alias}
\newcommand{\New}{C:/Users/yuzhe/Dropbox/bib/New}
\newcommand{\Main}{C:/Users/yuzhe/Dropbox/bib/Main}
\newcommand{\FP}{C:/Users/yuzhe/Dropbox/bib/FP}

\newcommand{\CP}{\mathcal{P}} %
\newcommand{\CG}{\mathcal{G}} %
\newcommand{\CE}{\mathcal{E}} %
\newcommand{\CV}{\mathcal{V}} %
\newcommand{\ctl}{{\rm ctl}} 

\newcommand{\TiaQED}{\hfill $\triangle$} 

\newcommand{\R}{\mathbb{R}} %

\newcommand{\diag}{{\rm diag}} %
\newcommand{\sign}{\operatorname{sign}}

\newcommand*{\QE}{\hfill\ensuremath{\square}}%
\newenvironment{pfof}[1]{\vspace{1ex}\noindent{\textit{Proof of
			#1:}}\hspace{0.5em}} {\hfill\QE\vspace{1ex}}

\usepackage{amsthm}

\usepackage{xcolor}         
\usepackage[colorlinks]{hyperref}
\hypersetup{citecolor=[RGB]{5,112,176},linkcolor=[RGB]{5,112,176},urlcolor=[RGB]{5,112,176}}

\theoremstyle{definition}
\newtheorem{theorem}{Theorem}

\newtheorem{lemma}{Lemma}

\theoremstyle{remark}
\newtheorem{remark}{Remark}
\newtheorem{example}{Example}
\newtheorem{assumption}{Assumption}
\newtheorem{definition}{Definition}

\newcommand{\orcidicon}[1]{%
    \href{https://orcid.org/#1}{\includegraphics[height=1.5ex]{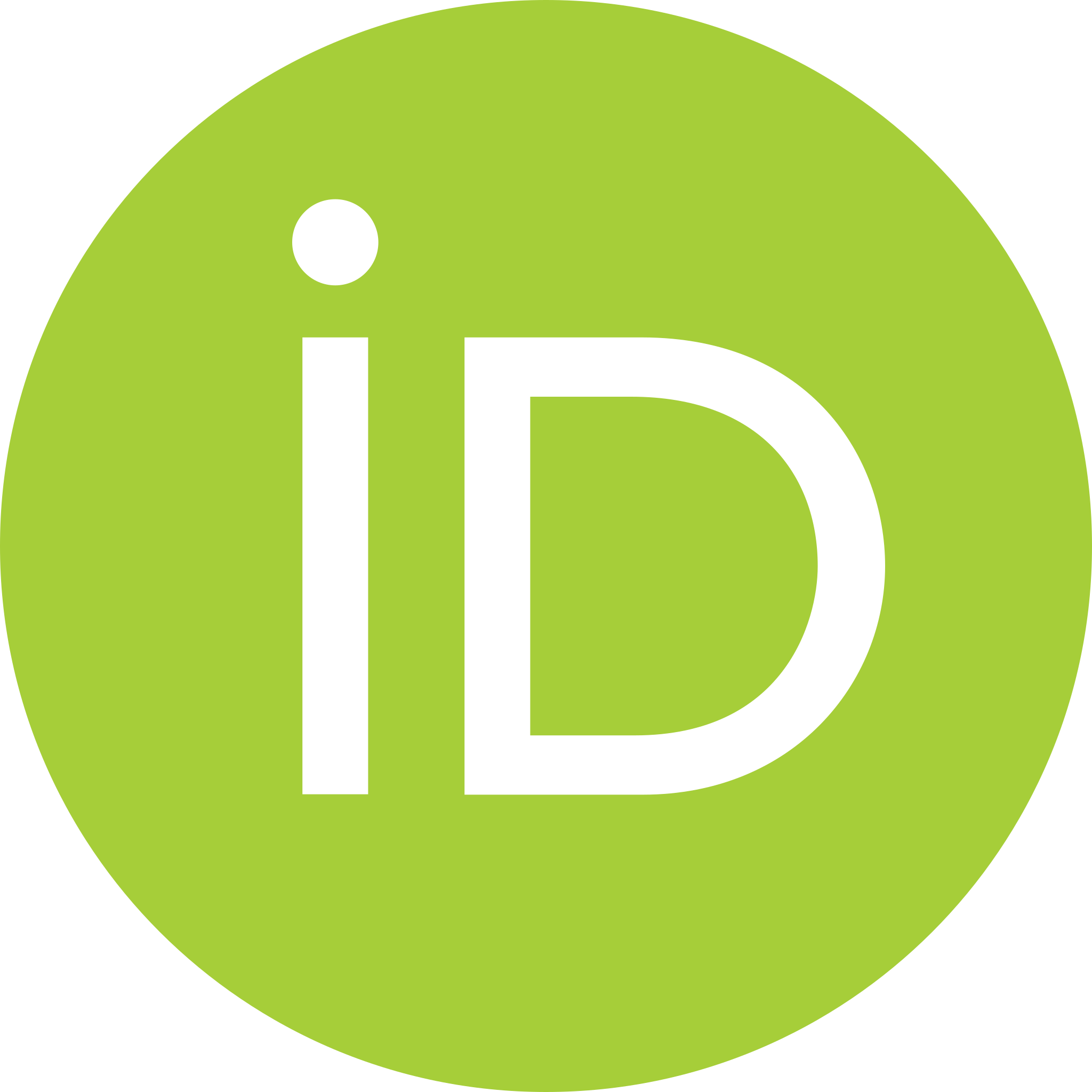}}%
}

\title{\Large \bf Vibrational Stabilization of Complex Network Systems}

\author{ Alberto Maria Nobili~\orcidicon{0000-0002-3073-8503}, 
	Yuzhen Qin~\orcidicon{0000-0003-1851-1370}, 
	Carlo Alberto Avizzano~\orcidicon{0000-0001-5802-541X}, 
	Danielle S. Bassett~\orcidicon{0000-0002-6183-4493}, \\
	and Fabio Pasqualetti~\orcidicon{0000-0002-8457-8656} 
  \thanks{This research was funded in part by Aligning Science Across
    Parkinson’s (ASAP-020616) through the Michael J. Fox Foundation for Parkinson’s Research (MJFF) and in part by NSF (NCS-FO-1926829) and
    ARO (W911NF1910360). 
    A. M. Nobili, Y. Qin and F. Pasqualetti are
    with the Department of Mechanical Engineering, University of
    California, Riverside (a.nobili5@studenti.unipi.it, \{yuzhenqin,
    fabiopas\}@engr.ucr.edu). 
    C. A. Avizzano is with the Perceptual Robotics Laboratory at the IIM Institute, Department of Excellence in Robotics and A.I., Scuola Superiore Sant’Anna, Pisa, Italy (carlo@sssup.it). 
    D. S. Bassett is with the Department of Bioengineering, the Department of Electrical \& Systems Engineering, the Department of Physics \& Astronomy, the
    Department of Psychiatry, and the Department of Neurology,
    University of Pennsylvania, and the Santa Fe Institute
    (dsb@seas.upenn.edu). 
    Y. Qin, F. Pasqualetti, and D. S. Bassett are also
    affiliated with Aligning Science Across Parkinson’s (ASAP)
    Collaborative Research Network, Chevy Chase, MD20815.
	The manuscript was prepared using LaTex on Texstudio, and the source code can be downloaded at: \hyperref[a]{https://arxiv.org/abs/2308.05823}. All simulations were performed in MATLAB, and the corresponding code is available at this repository: \hyperref[b]{https://doi.org/10.5281/zenodo.17179335}.} }

\begin{document}
	
\maketitle
\thispagestyle{empty}
\pagestyle{empty}

\begin{abstract}
	Many natural and man-made network systems need to maintain certain patterns, such as working at equilibria or limit cycles, to function properly. Thus, the ability to stabilize such patterns is crucial. Most of the existing studies on stabilization assume that network systems' states  can be measured online so that feedback control strategies can be used. However, in many real-world scenarios, systems' states, e.g., neuronal activity in the brain, are often difficult to measure. In this paper, we take this situation into account and study the stabilization problem of linear network systems with an open-loop control strategy---\textit{vibrational control}.  We derive a graph-theoretic sufficient condition for \textit{structural} vibrational stabilizability, under which network systems can always be stabilized. We further provide an approach to select the locations in the network for control placement and design corresponding vibrational inputs to stabilize systems that satisfy this condition. Finally, we provide some numerical results that demonstrate the validity of our theoretical findings. 
	
	
	
\end{abstract}

\section{Introduction}

Many natural and technological systems, such as gene regulation, neural circuits, and electric power grids, consist of large-scale interacting units. They are often modeled by complex network systems. Such network systems need to operate at certain equilibria or limit cycles  to function well. Therefore, guaranteeing their stability is vital. Loss of stability may lead to blackout in power grids \cite{JWB-KB-MDB-AQG-EG-VR-JDR-SS-CAS-MEW:21a} or neurological disorders in the brain \cite{JP-DCM-JJGR-SCA:2013}. For instance, the loss of the stability of normal coordinated brain activity leads to increased  synchrony in the basal ganglia and exaggerated phase-amplitude coupling in motor cortex, which are closely associated with Parkinson disease \cite{DHC-RWES:2013a}. Therefore, it is fundamental to be able to stabilize desired dynamic patterns of such network systems.  Most of existing studies on stabilization rely on the assumption that real-time states can be measured so that feedback control strategies can be used. However, in many real-world scenarios, states cannot be observed or measured directly. For instance, existing techniques find difficulty in precisely measuring neuronal activity in the brain, which poses challenges to restoring the stability of certain patterns of brain activity using feedback-based treatments. 


Vibrational control is a strategy to control a system without measuring its states \cite{SMM:80}. By injecting pre-designed high-frequency signals, it can stabilize  various engineering systems, e.g., inverted pendulums, chemical reactors, and under-actuated robots (see \cite{REB-JB-SMM:86b,BS-BTZ:97,CX-TY-MI:2018} and the references therein). It may also explain the mechanism of deep brain
stimulation \cite{YQ-DSB-FP:22a}, a neurosurgical technique used to treat several brain disorders including Parkinson disease. In this paper, we show how network systems can be stabilized by vibrational control. We concentrated on linear dynamics since stability of equilibria and limit cycles in nonlinear networks can often be studied by analyzing their linearized counterpart.

\textbf{Related work}. 
Stabilizability of linear network systems has attracted many interests  (e.g., see \cite{AC-MM:15, TJ-THL-18}). Recent works have studied structural stabilizability of network systems, where the network structure plays a central role in determining the stabilizability of a system \cite{SP-SK-PAA:15,LJ-CX-PS-PGJ-PVM:19}. Some studies have investigated the stabilizability of networks under malicious attacks (e.g., see \cite{DPC-TP:15, DA-SF-DMD:16}). Controllability of network systems has also received extensive attention in the past decades \cite{FP-SZ-FB:13q,SSM-MH-MM:17}; structural controllability is one of the most well-studied problems (e.g., \cite{JJ-HJW-HLT-KMC:19,GB-KO-AS-KM:21,AW-SM-YY-KX:22}). All the above studies assume that the states or outputs are measurable. However, this paper aims to stabilize network systems with an open-loop strategy, vibrational control, without that assumption. To the best of our knowledge, it is the first one to study vibrational stabilization of linear network systems. 






\textbf{Paper contribution}. The main contribution of this paper is threefold. First, we obtain a sufficient graph-theoretic condition for the vibrational stabilizability of linear network systems. Specifically, we find that for an arbitrarily parameterized system, if removing all the bidirected edges of the network associated with it results in a network that contains no cycles, this system is vibrationally stabilizable.  Second, we present a method to design vibrational control that targets a part of the edges in  the network to stabilize systems satisfying the aforementioned condition. Specifically, we propose an algorithm to place control inputs and we also show how to configure the  frequencies and amplitudes of the corresponding sinusoidal vibrations. Third, using averaging techniques, we define a notion of \textit{functional system} for the vibrationally controlled system. We further find that the working mechanism of vibrational stabilization in network systems is to \textit{functionally} modify the network parameters, such as changing the edge weights or removing edges. Finally, some numerical studies are also performed to validate our theoretical findings. 


\textbf{Notation}. Let $\R$ denote the set of real numbers. Given a directed graph $\CG=(\CV,\CE)$, denote the edge from $i$ to $j$ as $(i,j)$. We use $i_1\to i_2 \to \dots \to i_{k-1}\to i_k$ to denote a directed path from $i_1$ to $i_k$ passing through the nodes $i_2,\dots,i_{k-1}$.  Given $d_1,d_2,\dots,d_n\in \R$, $\diag(d_1,\dots,d_n)$ is the diagonal matrix. 


\section{Problem formulation}

\subsection{Linear network systems}
Consider a network represented by the directed graph $\CG :=(\CV,\CE)$, where $\CV=\{1,2,\dots,n\}$ and $\CE \subseteq \CV \times \CV$ are the sets of nodes and edges, respectively. Let $a_{ij}\in \R$ be the weight of the edge $(j,i)\in \CE$, and define the \textit{weighted adjacency} matrix as $A=[a_{ij}]_{n\times n}$, where $a_{ij}=0$ whenever $(j,i)\notin \CE$. Now, consider the linear network system described by
\begin{equation}\label{original_networked_sys}
	\dot{x}(t) = (D +A)  x(t):=Mx(t), 
\end{equation}
where $x_i\in\R$ in $x:=[x_1,x_2,\dots,x_n]^\top$ and $d_i \in \R$ in $D:=\diag(d_1,d_2,\dots,d_n)$ are the state and the intrinsic dynamics of the node $i$, respectively. In this paper, we are interested in the situation where $\sign (a_{ij})=\sign (a_{ji})$ whenever $a_{ij}\neq 0$ and $a_{ji}\neq 0$. Also, we consider that the intrinsic dynamics of each node is stable, i.e., the following assumption holds.

\begin{assumption}\label{assumption:stable}
	Assume that  $d_i<0$ for all $i=1,2,\dots,n$.
\end{assumption}

Yet, this assumption does not ensure that the overall system \eqref{original_networked_sys} is also stable. The network system can become unstable just because of connections. We illustrate this point in the following example.


\begin{example}
	Consider a network system associated with the graph depicted in Fig.~\ref{fig:ex_destabilization}~(a). As in Assumption~\ref{assumption:stable}, individual node systems are set to be stable. When the coupling strength $k$ is small, it can be observed from Fig.~\ref{fig:ex_destabilization}~(b) that the overall system is stable; when $k$ becomes larger, the system becomes unstable (see Fig.~\ref{fig:ex_destabilization}~(c)).  \TiaQED
\end{example}

There are certainly other scenarios than the one in this example where networks of stable units become unstable. 
In this paper, we aim to investigate whether and how network systems described by \eqref{original_networked_sys} can be stabilized by a classic open-loop control strategy: \textit{vibrational control}.

\subsection{Vibrational control in network systems}

Given a linear system $\dot x =M x$, vibrational control introduces vibrations to the system matrix $M$, resulting in the following controlled system
\begin{equation}\label{controlled_net_compact}
	\dot x = \Big(M+\frac{1}{\varepsilon}V\big(\frac{t}{\varepsilon}\big) \Big) x, 
\end{equation} 
where the \textit{zero-mean} control input $V(t)$ is often chosen to be periodic or quasiperiodic \cite{SMM:80,REB-JB-SMM:86a,REB-JB-SMM:86b}. For instance, a widely-used $V(t)=[v_{ij}(t)]_{n\times n}$ has $v_{ij}(t)=\mu_{ij}\sin(\omega_{ij} t)$ for some constant $\mu_{ij}\approx 1$ and $\omega_{ij}\approx 1$. The parameter $\varepsilon>0$ determines the frequency of the vibrations. An appropriate configuration of vibrations can stabilize an unstable system without any measurements of the states \cite{SMM:80}. 

\begin{figure}[t]	
	\centering
	\includegraphics[width=0.42\textwidth]{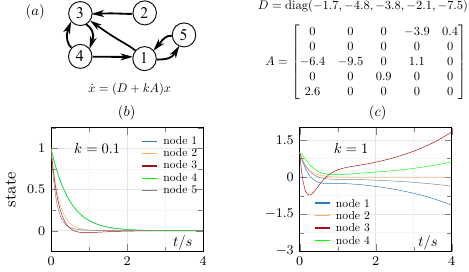}	
	\caption{Coupling strengths can shape the stability of network systems. (a) The network structure and the system dynamics. (b) If the coupling strength is small (e.g., $k=0.1$) the system is stable. (c) If the coupling strength is large (e.g., $k=1$), the system becomes unstable.}
	\label{fig:ex_destabilization}
	\vspace{-10pt}
\end{figure}

For general linear systems, one can introduce vibrations to any $m_{ij}$ in the system matrix $M=[m_{ij}]_{n\times n}$. When it comes to network systems, one can no longer do the same. For instance, one cannot introduce vibrations to $m_{12}$ for the network system depicted in Fig.~\ref{fig:ex_destabilization}. This is because there is no connection between the nodes $1$ and $2$, and it is unreasonable to inject vibrations in a nonexistent  connection.

Therefore, vibrational control in network systems needs to be constrained by the network structure. It is then natural to assume that vibrations can only be introduced to the intrinsic dynamics of node systems $d_i$ and the coupling strengths $a_{ij}$ in \eqref{original_networked_sys}. As a result, the vibrational control matrix $V(t)$ has the following constraint:
\begin{align}\label{constrain}
	&v_{ij}(t)=0, &\text{if } m_{ij}=0.
\end{align}
In other words, the vibrational control needs to have the same \textit{sparsity pattern} of the matrix $M$.

Our goals in this paper becomes to: 1)  investigate the conditions under which a network linear system is stabilizable using the vibrational control with the above constraint, and, subsequently, 2) study how to design vibrational control to stabilize a system  satisfying such conditions.  

Following \cite{SMM:80}, we now generalize the definition of vibrational stabilizability to network systems.
\begin{definition}
	The network system described by \eqref{original_networked_sys} is \textit{vibrationally stabilizable} if there exists a vibrational input $V(t)$ that satisfies \eqref{constrain} and such that the controlled system \eqref{controlled_net_compact} is asymptotically stable.
\end{definition}

\section{Vibrational Control of Network Systems}\label{sec:control}

\subsection{Averaged system and functional network}
To study vibational control, a key step is to analyze the stability of the controlled system \eqref{controlled_net_compact}. Since $\eqref{controlled_net_compact}$ is a time-varying system, a typical approach is to associate it with an averaged system. Then, the stability of \eqref{controlled_net_compact} can be indirectly studied by investigating its time-invariant averaged counterpart (e.g., see \cite{SMM:80, CX-TY-MI:2018}).\\

The first step is to change the timescale to $s=t/\varepsilon$, so that the system \eqref{controlled_net_compact} becomes
\begin{equation}\label{controlled_net_compact_s}
	\frac{dx}{ds} = \Big(\varepsilon M+V(s) \Big) x, 
\end{equation} 

Now, the standard first-order averaging (e.g., see \cite[Chap. 10]{HKK:02-bis}) is not applicable here. Indeed, since $V(s)$ has zero mean, applying the first-order averaging to \eqref{controlled_net_compact_s} just eliminates the $V(s)$ term and results in the uncontrolled system ${dx/ds = \varepsilon Mx}$.

To avoid this issue, we change the coordinates of \eqref{controlled_net_compact} before using averaging. Specifically, we introduce an auxiliary system
\begin{equation*}
	\frac{dx}{ds}=V(s)x(s)
\end{equation*}
and let $\Psi(s,s_0)$ be its state transition matrix. Applying the change of coordinates $z(s)=\Psi(s,s_0)^{-1}x(s)$, the system \eqref{controlled_net_compact_s} can be rewritten as
\begin{equation} \label{controlled_net_compact_z}
	\frac{dz}{ds}=\varepsilon \Psi^{-1}(s,s_0) M \Psi(s,s_0) z.
\end{equation}

Since $V(s)$ is often (quasi-)periodic, $\Psi$ is bounded. Then, the stability of \eqref{controlled_net_compact_z} implies that of \eqref{controlled_net_compact}.

We then introduce the averaged system of \eqref{controlled_net_compact_z}:
$
	\frac{dz}{ds}= \varepsilon \bar{M}z,
$
where 
$
	\bar{M} = \lim_{T\rightarrow\infty}\frac{1}{T}\int_{0}^{T} \Psi(s,s_0)^{-1}M\Psi(s,s_0)ds.
$

Changing  the timescale back to $t=\varepsilon s$ leads to
\begin{equation} \label{controlled_net_compact_z_avg}
	\dot z= \bar{M}z.
\end{equation}

The following lemma provides the relation between the stability of the original and averaged systems  \eqref{controlled_net_compact_s} and \eqref{controlled_net_compact_z_avg}. The proof follows the same line as in \cite{SMM:80,CX-TY-MI:2018,QY-KY-BDOA-CM:2021}. 
\begin{lemma}\label{averge_to_original}
	Assume that $\bar M$ of the averaged system \eqref{controlled_net_compact_z_avg} is Hurwitz.  Then, there exists $\varepsilon^*>0$ such that for any $\varepsilon> \varepsilon^*$, the system \eqref{controlled_net_compact_s} is asymptotically stable. \TiaQED
\end{lemma}

This lemma implies that a system is stabilizable if there exists a vibrational control such that the averaged system is stable. Also, to stabilize a system, the problem reduces to find a vibrational configuration $V(s)$ such that the averaged system \eqref{controlled_net_compact_z_avg} is stable. Therefore, we will refer to the system \eqref{controlled_net_compact_z_avg} as the \textit{functional system} of the controlled system \eqref{controlled_net_compact_s}. 

Now, let us rewrite the functional system \eqref{controlled_net_compact_z_avg} into
\begin{equation} \label{averaged:separate}
	\dot z = \bar{M}z=(\bar D +\bar A)z,
\end{equation}
where $\bar D=\diag(\bar d_1,\dots,\bar d_n)$ is the diagonal matrix of $\bar M$ and $\bar A=\bar M -\bar D$ is the off-diagonal matrix. 

This functional system can be also taken as a network system. Its differences from the original network system \eqref{original_networked_sys} are: 1) the directed graph  associated with \eqref{averaged:separate} becomes $\bar \CG:=(\CV,\bar \CE)$, where the weighted adjacency matrix is described by the matrix $\bar A$, and 2) the intrinsic dynamics of node systems become $\bar D$.  We refer to the network described by $\bar \CG:=(\CV,\bar \CE)$ as the \textit{functional network} of the controlled system. As one may have observed, vibrational control can introduce the following functional changes to a network system: 1) modification of the intrinsic dynamics, and 2) alteration of the network weights or structure.


\subsection{Vibrational stabilizability}
\begin{figure}[t]
	\centering
	\includegraphics[scale=1.1]{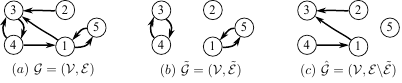}
	\caption{Illustration of the definitions. (a) The original graph $\CG$. (b) The graph composed of bidirected edges.  (c) The unidirected residual graph $\hat \CG$ is a DAG. From Theorem~\ref{th:stabilizability}, any linear network system associated with $\CG$ is vibrationally stabilizable.}
	\label{fig:application_ex}
\end{figure}

We present our main results in this subsection. First, we provide some relevant definitions (see Fig.~\ref{fig:application_ex} for an illustration).

\begin{definition}\label{defi:DAG}
	A {directed acyclic graph (DAG)} is a directed graph that does not contain any directed cycle.
\end{definition}

\begin{definition}\label{def:bidirected}
	Any two nodes $i,j$ in $\CG=(\CV,\CE)$ are said to be connected bidirectionally if $(i,j)\in \CE$ and $(j,i)\in \CE$. The two edges $(i,j)$ and $(j,i)$ are referred to as \textit{bidirected edges}. Let $\tilde{\CE}:=\{(i,j)\in\CE:(j,i)\in \CE\}$ be the set of all bidirected edges. 
\end{definition}

\begin{definition}\label{residual}
	The graph $\hat \CG=(\CV,\hat \CE)$ with $\hat \CE=\CE\backslash \tilde \CE$ is said to be the \textit{unidirected residual} of $\CG=(\CV,\CE)$.
\end{definition}


\begin{theorem}[\textbf{Structural vibrational stabilizability}]\label{th:stabilizability}
	Consider a network system described by \eqref{original_networked_sys} that satisfies Assumption \ref{assumption:stable} and is associated with the graph $\mathcal{G}=(\mathcal{V}, \mathcal{E})$. It is vibrationally stabilizable if the unidirected residual $\hat \CG$ of $\mathcal{G}$ is a DAG. \TiaQED
\end{theorem}

Note that the condition  is graph-theoretic, only depending on the network structure, and the weights of the edges do not matter. This is why we call it a structural vibrational stabilizability condition. The following example illustrate how Theorem~\ref{th:stabilizability} can be applied.

\begin{example}\label{example:stabilizable}
	Consider a system that is associated with the network depicted in Fig.~\ref{fig:application_ex}~(a). Since the unidirected residual $\hat \CG$ contains no directed cycles (see Fig.~\ref{fig:application_ex}~(c)), this system is vibrationally stabilizable. We stress that the condition for vibrational stabilizability in Theorem~\ref{fig:application_ex} is a sufficient one. A network system that does not satisfy it may still be stabilized. We will discuss this point in Example~\ref{exam:not_necessary}. \TiaQED
\end{example}



As we mentioned earlier, a system is stabilizable if one can find a vibrational control that actually stabilizes it.
Next, we consider a particular form of vibrational control to facilitate the proof of Theorem~\ref{th:stabilizability} and to show how to design vibrational control. 

\subsection{Design of vibrational control}
We consider sinusoidal vibrations, which means that $V(s)=[v_{ij}(s)]_{n \times n}$ in the system \eqref{controlled_net_compact_s} has the following form:
\begin{equation}\label{sinu:v_ij} 
	v_{ij}(s) = \mu_{ij}\sin(\omega_{ij}s). 
\end{equation}
Also, in this paper, we consider that vibrations are only introduced to the edges in the network $\CG=(\CV,\CE)$. Let $\CE_{\ctl}\subseteq \CE$ be the target set of edges that control inputs are injected to. Then, the vibrations in \eqref{sinu:v_ij}  satisfy
\begin{equation}\label{individual:v_ij} 
	\mu_{ij}
	\begin{cases}
		\neq 0,\qquad \text{if } (j,i)\in \CE_\ctl,\\
		=0, \hfill \text{otherwise}.
	\end{cases}
\end{equation}

The following theorem provides an approach to design the vibrational control to stabilize a system (which also implies that the proof of Theorem~\ref{th:stabilizability} follows directly). Without loss of generality, we assume that there are $m$ pairs of bidirected edges in $\CG$, each denoted as $\CP_k:=\{(i_k,j_k),(j_k,i_k)\}$, $k=1,\dots,m$.

\begin{algorithm}[t]
	\caption{Control Input Placement}
	\label{control:placement}
	\begin{algorithmic}[1]
		\State \textbf{Input:}  $\CG$, unidirected residual $\hat \CG$, and $\CP_k,k=1,\dots,m$
		\State \textbf{Initialize:} $\tilde{\CE}_2=\{\}$, $\bar \CG=\hat \CG$, and $\bar \CE=\hat \CE$
		\For{$k=1:m$}
		\State	find $e_k\in \CP_k$ such that $\bar\CG=(\CV,\bar \CE \cup \{e_k\} )$ is DAG
		\State	update $\bar \CE =\bar \CE\cup \{e_k\}$, $\bar \CG=(\CV,\bar \CE)$, $\tilde{\CE}_2=\tilde{\CE}_2\cup\{e_k\}$
		\EndFor
		\State $\tilde{\CE}_1=\tilde{\CE}\backslash \tilde{\CE}_2$
	\end{algorithmic}
\end{algorithm}

\begin{theorem}[\textbf{Design of vibrational control}]\label{The:control}
	Consider a network system described by \eqref{original_networked_sys} that is associated with the graph $\mathcal{G}=(\mathcal{V}, \mathcal{E})$, and assume it satisfies the conditions in Theorem~\ref{th:stabilizability}. Let the control target set be $\CE_{\ctl}=\tilde{\CE}_1$, where $\tilde{\CE}_1$ is generated by Algorithm~\ref{control:placement}. Then, there exist constants $\mu_{ij},\omega_{ij}$ such that the vibrational control defined in \eqref{sinu:v_ij} and  \eqref{individual:v_ij} leads to the following two statements:
	
	(I) The functional dynamics of the controlled system, described by \eqref{averaged:separate}, is asymptotically stable. Further, the functional intrinsic dynamics  satisfy $\bar D=D$, and the functional network is represented by a DAG $\bar \CG=(\CV,\CE\backslash \CE_{\ctl})$.
	
	(II) There exists $\varepsilon^*>0$  such that the controlled network system \eqref{controlled_net_compact} is asymptotically stable for any $\varepsilon>\varepsilon^*$.\TiaQED
\end{theorem}

The statement (I) of Theorem~\ref{The:control} states that the vibrations preserve the intrinsic dynamics of node systems, but functionally remove the edges in the target set $\CE_{\ctl}$ from the original network, resulting in a \textit{directed acyclic} functional network. The following lemma guarantees that such a functional system is asymptotically stable (the proof can be found in the Appendix). 

\begin{lemma}\label{lemma:stable_graph}
	Assume that the network system described in  \eqref{original_networked_sys} satisfies Assumption \ref{assumption:stable} and is associated with the graph $\mathcal{G}=(\mathcal{V}, \mathcal{E})$. Then, this system is asymptotically stable if the graph $\mathcal{G}$ is a DAG.
\end{lemma}

Algorithm~\ref{control:placement} provides an approach to select the control target set $\CE_{\ctl}$ that contains the edges that we want to functionally remove. Given $m$ pairs of bidirected edges, our goal is to remove one edge from each pair such that the remaining graph becomes a DAG. To decide which edges to remove (i.e., $\tilde{\CE}_1$), we first decide which edges to keep (i.e., $\tilde{\CE}_2$). The key idea is to add back $m$ directed edges to the unidirected residual graph $\hat \CG$ in $m$ steps, one from each pair at each step. The following lemma ensures that, starting from a directed acyclic $\hat \CG$, the graph with one new edge added at each step in Algorithm~\ref{control:placement} is always a DAG. After $m$ steps, the resulting graph $\bar \CG$ is also a DAG. Then, it subsequently becomes clear which set of edges to remove.

\begin{lemma}\label{acyclic:adding}
	Consider a DAG $\CG=(\CV,\CE)$. For any pair of nodes $i,j\in \CV$ satisfying $(i,j)\notin \CE$ and $(j,i)\notin \CE$, there exists a directed edge $e\in \{(i,j),(j,i)\}$ such that the graph $\CG'=(\CV,\CE \cup \{e\})$ is still a DAG. \TiaQED
\end{lemma}

\begin{figure}
	\centering
	\includegraphics[scale=0.5]{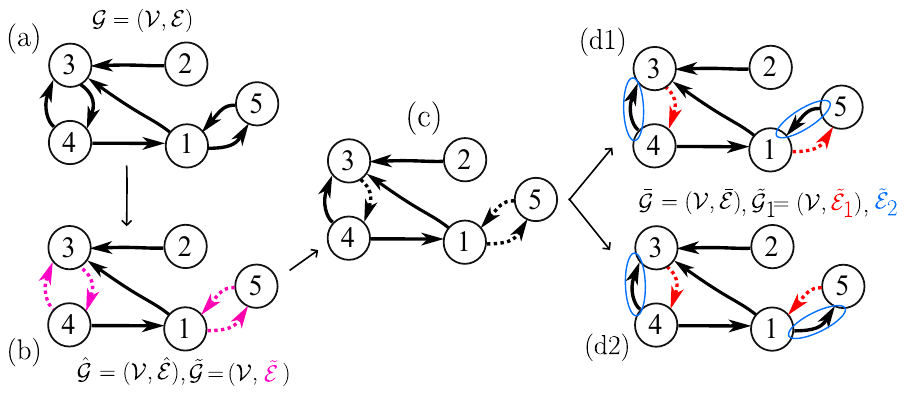} 
	\caption{Illustration of how to select a target set of edges to  control. (a) Original graph $\CG$. (b) The set of bidirected edges $\tilde{\CE}$ (dashed) and the unidirect residual graph $\hat \CG$. From (b) to (c), we add one directed edge between nodes $3$ and $4$, resulting in a DAG; from (c) to (d), we add another edge between nodes $1$ and $5$, still ensuring the graph is a DAG. In (d1) and (d2), we show that there is more than one way to add such edges.  The added edges form the set $\tilde{\CE}_2$, and then the control target set becomes $\tilde{\CE}_1=\tilde{\CE}\backslash \tilde{\CE}_2$. Note that both the graphs $\bar \CG$ and $\tilde{\CG}_1$ are DAGs. }	
	\label{illustration:algo}
	\vspace{-10pt}
\end{figure}

\begin{remark}
	We provide an example in Fig.~\ref{illustration:algo} to illustrate how a control target set $\tilde{\CE}_1$ can be selected. We emphasize that, for a given network, there can be multiple choices of $\tilde{\CE}_1$, which means that there is more than one way to inject the vibrational control (see Fig.~\ref{illustration:algo} (d1) and (d2) for an example). 
\end{remark}

Next, we provide the proof of Theorem~\ref{The:control}, which needs the following lemma (the proof can be found in the Appendix).

\begin{lemma}\label{lemma:elimination_of_edges}
	For a network system described by \eqref{original_networked_sys} with the associated graph $\mathcal{G}=(\mathcal{V}, \mathcal{E})$, let $\tilde \CE$ be as in Definition~\ref{def:bidirected}. Consider any subset $\tilde \CE_1\subseteq\tilde \CE$ and let  $\bar \CE= ({\CE}\backslash\tilde{\CE}_1)$. If $\tilde{\CG}:=(\CV,\tilde \CE_1)$ is a DAG, then there exist constants $\mu_{ij},\omega_{ij}$ such that the vibrational control defined in \eqref{sinu:v_ij} and  \eqref{individual:v_ij} leads to a functional network represented by $\bar \CG=(\mathcal{V}, \bar \CE)$. \TiaQED
\end{lemma}

\begin{pfof}{Theorem~\ref{The:control}}
	From Lemma~\ref{acyclic:adding}, the graph $\bar \CG=(\CV,\hat \CE \cup \tilde{\CE}_2)$ is directed acyclic.  Then, it holds that the edges in $\tilde{\CE}_2$ do not form a directed cycle. Consequently, since the set $\tilde{\CE}_1$ satisfies $\tilde{\CE}_1=\tilde \CE\backslash \tilde{\CE}_2$, $\tilde{\CE}_1$ also contains no cycles. It follows from Lemma~\ref{lemma:elimination_of_edges} that there exist constants $\mu_{ij},\omega_{ij}$ such that the vibrational control defined in \eqref{sinu:v_ij} and  \eqref{individual:v_ij}  leads to a functional dynamics with function network represented by $\bar \CG=(\mathcal{V}, \bar \CE)$. Lemma~\ref{lemma:stable_graph} implies the asymptotic stability of the functional system, which provide the statement (I). The statement (II) follows directly from Lemma~\ref{averge_to_original}.	
\end{pfof}

\subsection{Numerical Studies}
In this subsection, we use two numerical examples to demonstrate our theoretical results.

\begin{figure}[t]
	\centering
	\includegraphics[scale=0.8]{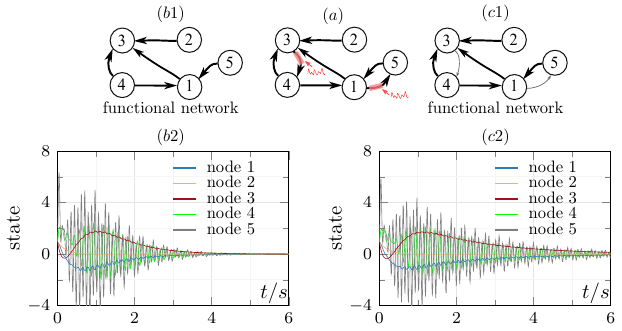}
	\caption{Vibrational stabilization of the  network system in Fig.~\ref{fig:ex_destabilization} (c). (a) The network associated with the uncontrolled system (the red signals are the vibrations we want to inject). (b) By removing connections from the functional network [see (b1)], the vibrational control stabilizes the system [see (b2)]. (c) By weakening connections in the functional network [see (c1)], the system is also stabilized [see (c2)]. Vibrations are introduced to edges $\lambda_{43}$ and $\lambda_{51}$ in both (b) and (c).}
	\label{fig:ctr}
	\vspace{-2pt}
\end{figure}

\begin{example} \label{exam:sufficiency}
	First, we revisit the unstable network system in  Fig.~\ref{fig:ex_destabilization} (i.e., $k=1$). As argued in Example~\ref{example:stabilizable}, this system is vibrationally stabilizable. Following Theorem~\ref{The:control} and Fig.~\ref{illustration:algo}, we inject vibrations to the edges $(3,4)$ and $(1,5)$. Specifically, $v_{43}(t) = 50\sqrt{1.8/1.1}\sin(50t)$, $v_{51}(t) = 100\sqrt{6.5} \sin(50\sqrt{2}t)$, and $v_{ij}(t)=0$ for any other $i,j$. As shown in Fig.~\ref{fig:ctr}~(b), these vibrations have functionally removed the edges $(3,4)$ and $(1,5)$, resulting in a stable system. Furthermore, if we decrease the amplitudes of the vibrations,  the system can be stabilized without removing the edges in the functional network (see Fig.~\ref{fig:ctr}~(c)).

\end{example}

\begin{figure}[t]
	\centering
	\includegraphics[scale=0.9]{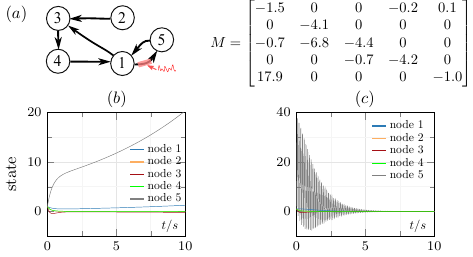}
	\caption{Example of a vibrational stabilizable network system that does not satisfy Theorem~\ref{th:stabilizability}. (a) The network structure and the system matrix. (b) Without vibrational control, the system is unstable. (c) With vibrational control on edge $\lambda_{51}$, the system becomes stable even if there is a cycle involving nodes $1, 3$ and $4$.}
	\label{fig:ctr_exc}
	\vspace{-18pt}
\end{figure}

\begin{example} \label{exam:not_necessary}
	Now, we consider another unstable network system depicted in Fig.~\ref{fig:ctr_exc}~(a). One can observe that removing the bidirected edges leads to a network that contains a directed cycle. Thus, the conditions in Theorem~\ref{th:stabilizability} are violated. However, a carefully designed vibration injected to the edge $(5,1)$ stabilizes the system (see Fig.~\ref{fig:ctr_exc}~(c)). Specifically, $v_{51}(t) = 100\sqrt{2.6/0.4} \sin(50\sqrt2 t)$, and $v_{ij}(t)=0$ for any other $i,j$. With this example, we wish to mention that the conditions in Theorem~\ref{th:stabilizability} are just sufficient. It remains of interest to investigate the necessary and sufficient conditions.
	
\end{example}

\subsection{Vibrational control can improve robustness}
In the previous subsections, we have studied how vibrational control can be used to stabilize unstable network systems. Now, we show that it can also improve robustness of stable network systems.

Following \cite{DH-AJP:86}, we employ the Unstructured Real Stability Radius (URSR) to measure the robustness of the linear network system \eqref{original_networked_sys}. Specifically, the URSR of the system \eqref{original_networked_sys} is defined as
\begin{equation}
	r_\R(M) := \inf_{\Delta\in\R^{n\times n}}\{||\Delta||_2 : \alpha(M+\Delta)\ge0\},
\end{equation}
where $\alpha(\cdot)$ denotes the spectral abscissa\footnote{The spectral abscissa of a square matrix is the largest real part of its eigenvalues.} of a matrix, and $||\cdot||_2$ is the spectral norm. The URSR provides a worst-case measure for the robustness of a system in the sense that all perturbations with $||\Delta||<r_\R(M)$ are guaranteed to preserve the stability of the perturbed system.

\begin{figure}[t]
	\centering
	\includegraphics[scale=0.73]{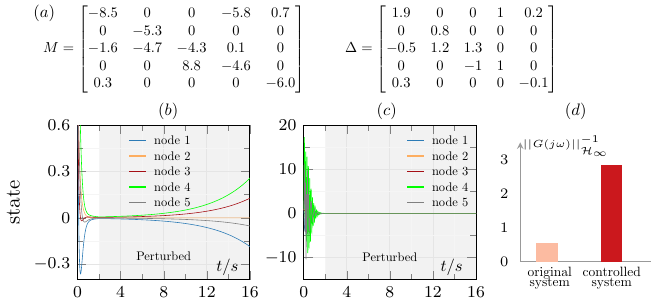}
	\caption{Vibrational control improves robustness. (a) The system matrix $M$, and the perturbation matrix $\Delta$. (b) Without vibrational control, the system becomes unstable in presence of perturbation. (c) With vibrational control on edge $\lambda_{51}$ and $\lambda_{43}$, the system remains stable even in presence of perturbation. Perturbation are introduced at $t=2 s$ (d) A bar chart showing lower bounds for the URSRs.}
	\label{fig:robustness}
	\vspace{-18pt}
\end{figure}

We employ the $\mathcal{H}_{\infty}$ norm to roughly approximate URSR, with the relation between them \cite{LQ-BB-AR-EJD-PMY-JCD:95} given by
\begin{equation}
	r_\R(M)\ge\left[\sup_{\omega\in\R}\bar{\sigma}\left((j\omega I-M)^{-1}\right)\right]^{-1}
	=\rVert G_M(j\omega) \rVert_{\mathcal{H}_\infty}^{-1},
\end{equation}
where $\bar{\sigma}(\cdot)$ denotes the maximum singular value, $j=\sqrt{-1}$ is the imaginary unit, and $G_M(j\omega)=(j\omega I-M)^{-1}$. We consider a stable system with the associated system matrix given in Fig.~\ref{fig:robustness}~(a).  To show how vibrational control can improve robustness, we introduce a perturbation on the intrinsic dynamics and the connections of the network system at time $t=2s$.  We compare the performance of the uncontrolled system and the controlled one (the vibrations introduced are the same as in Example \ref{exam:sufficiency}): the former loses  its stability due to perturbation, while the latter preserves it. 
	We also compare the inverse of the $\mathcal H_\infty$ norm of the original system and the functional one. As shown in Fig.~\ref{fig:robustness}~(d), the robustness of the latter is indeed improved by the vibrational control. 
	
	
	
	\section{Conclusions}
	In this paper, we study the vibrational stabilization of linear network systems. Different from vibrational control of general linear systems, vibrations are constrained by the network structure in our case. Sufficient conditions on the network structure are obtained such that any system associated with such networks are vibrationally stabilizable. We also provide an approach to design vibrational inputs to stabilize such systems. We put forth that the working principle of vibrational control is to functionally remove connections or modify the connection weights between node edges. We also present some numerical experiments to validate our theoretical findings. 
	As for future work, we are currently working to extend the results achieved to more general non-linear networks. 
	
	
	\appendix
	Here we present the proofs of the lemmas in Section~\ref{sec:control}. 
	
	%
	
	
	\begin{pfof}{Lemma~\ref{lemma:stable_graph}}
		As $\CG$ is also  a DAG,   according to \cite{JBJ-GG:00},  it can be topologically ordered. Therefore, one can  arrange the vertices of $\CG$ as a linear ordering that is consistent with all edge directions. In other words, there exists a permutation matrix $P$ such that the matrix $\hat M=:PMP^{-1}$ in the following system
		\begin{equation}\label{sys:permuted}
			\dot y =\hat M y , \hspace{1cm} \text{        with } y= Px,
		\end{equation}
		is lower-triangular. Under Assumption~\ref{assumption:stable}, one can derive that the diagonal entries of $\hat M$ are all negative, which means that $\hat M$ is Hurwitz. Therefore, the system \eqref{sys:permuted} is asymptotically stable, and so is the system \eqref{original_networked_sys}.		
	\end{pfof}
	
	\begin{pfof}{Lemma~\ref{acyclic:adding}}
		We construct the proof by contradiction. Now we assume that $\CG$ is directed acyclic and both $\CG'_1=(\CV,\CE \cup \{(i,j)\})$ and $\CG'_2=(\CV,\CE \cup \{(j,i)\})$ contain a cycle. For $\CG'_1$, there is a cyclic path 
		\begin{equation}\label{path:1}
			i_1 \to i_2\to \dots \to i \to j \to i_k \to \dots \to i_1.
		\end{equation}
		Likewise, for $\CG'_2$, there is also a cyclic path 
		\begin{equation}\label{path:2}
			j_1 \to j_2\to \dots \to j \to i \to j_k \to \dots \to j_1.
		\end{equation}
		One can observe from \eqref{path:1} and \eqref{path:2} that in the original graph $\CG$, there is a directed path from $i_1$ to $j_1$ and also from $j_1$ to $i_1$, which implies that $\CG$ contains a cycle. Observing that this is a contradiction completes the proof.
	\end{pfof}
	
	\begin{pfof}{Lemma~\ref{lemma:elimination_of_edges}}	
		As the graph $\tilde \CG_1 =(\CV, \tilde{\CE}_1)$ does not contain a directed cycle, there exists a permutation matrix $P$ such that $\tilde V(t):=P V(t) P^{-1}$ is lower-triangular (following the same line as in the proof of Lemma~\ref{lemma:stable_graph}). Considering the change of coordinates $y=Px$, then the controlled system becomes
		\begin{equation}\label{dynamics:y}
			\dot y = (\tilde M +\tilde{V}(t))y,
		\end{equation}
		where $\tilde M= P M P^{-1}$. 
		
		Now, observe that \eqref{dynamics:y} is a system controlled by a vibrational control that has a lower-triangular form. Then, letting $w_{ij}$ be incommensurable for different pairs of $i$ and $j$ and following the same steps as in \cite{SMM:80}, one can derive that the averaged system of \eqref{dynamics:y} is 
		\begin{equation*}
			\dot {y} = (\tilde M +B)y:= \bar M' y,
		\end{equation*}
		where $B= -\tilde M ^\top \odot C$ with $C=[c_{k \ell}]_{n\times n}$ such that $c_{kl}\ge 0$. Here, the value of each $c_{k \ell }$ is determined by the amplitude and frequency of the the vibrations (i.e., $\mu_{ij}$ and $\omega_{ij}$ in \eqref{sinu:v_ij}). 	Further, the definition of $\tilde{\CE}_1$ ensures that for any two nodes $i$ and $j$ such that $(j,i)\in \tilde{\CE}_1$, it holds that $(i,j)\notin \tilde{\CE}_1$ and $(i,j)\in {\CE}$, and thus $m_{ji}\neq 0$. This implies that for any $i$ and $j$ such that $\tilde{V}_{ij}(t)\neq 0$, it holds that $\tilde{m}_{ji}\neq 0$. Then, since $\sign (\tilde m_{ij})=\sign (\tilde m_{ji})$, one can choose a configuration of the amplitudes and frequencies, $\mu_{ij}$ $\omega_{ij}$, such that $\bar M'$ satisfies $\bar m'_{k\ell}=0$ for any vibrationally controlled $(\ell,k)$. 
		
		Subsequently, one can derive that the averaged system $\dot x=\bar M x$ of $\dot x =(M+V(t))x$ satisfies $\bar M = P^{-1}\tilde M' P$ and $\bar m_{ij}=0$  for any $(j,i)\in \tilde{\CE}_1$. As a consequence, the graph associated with this averaged system is $\bar \CG=(\CV,\CE\backslash \tilde{\CE}_1)$, which completes the proof.
		%
		%
		%
	\end{pfof}

	\bibliographystyle{IEEEtran}
	\bibliography{\alias,\FP,\Main,\New}

\begin{thebibliography}{10}
\providecommand{\url}[1]{#1}
\csname url@samestyle\endcsname
\providecommand{\newblock}{\relax}
\providecommand{\bibinfo}[2]{#2}
\providecommand{\BIBentrySTDinterwordspacing}{\spaceskip=0pt\relax}
\providecommand{\BIBentryALTinterwordstretchfactor}{4}
\providecommand{\BIBentryALTinterwordspacing}{\spaceskip=\fontdimen2\font plus
\BIBentryALTinterwordstretchfactor\fontdimen3\font minus
  \fontdimen4\font\relax}
\providecommand{\BIBforeignlanguage}[2]{{%
\expandafter\ifx\csname l@#1\endcsname\relax
\typeout{** WARNING: IEEEtran.bst: No hyphenation pattern has been}%
\typeout{** loaded for the language `#1'. Using the pattern for}%
\typeout{** the default language instead.}%
\else
\language=\csname l@#1\endcsname
\fi
#2}}
\providecommand{\BIBdecl}{\relax}
\BIBdecl

\bibitem{JWB-KB-MDB-AQG-EG-VR-JDR-SS-CAS-MEW:21a}
J.~W. Busby, K.~Baker \emph{et~al.}, ``Cascading risks: Understanding the 2021
  winter blackout in {T}exas,'' \emph{Energy Research \& Social Science},
  vol.~77, p. 102106, 2021.

\bibitem{JP-DCM-JJGR-SCA:2013}
P.~Jiruska, M.~De~Curtis \emph{et~al.}, ``Synchronization and desynchronization
  in epilepsy: Controversies and hypotheses,'' \emph{The Journal of
  Physiology}, vol. 591, no.~4, pp. 787--797, 2013.

\bibitem{DHC-RWES:2013a}
C.~De~Hemptinne, E.~S. Ryapolova-Webb \emph{et~al.}, ``Exaggerated
  phase--amplitude coupling in the primary motor cortex in parkinson disease,''
  \emph{Proceedings of the National Academy of Sciences}, vol. 110, no.~12, pp.
  4780--4785, 2013.

\bibitem{SMM:80}
S.~M. Meerkov, ``Principle of vibrational control: Theory and applications,''
  \emph{IEEE Transactions on Automatic Control}, vol.~25, no.~4, 1980.

\bibitem{REB-JB-SMM:86b}
R.~E. Bellman, J.~Bentsman, and S.~M. Meerkov, ``Vibrational control of
  nonlinear systems: Vibrational controllability and transient behavior,''
  \emph{IEEE Transactions on Automatic Control}, vol.~31, no.~8, pp. 717--724,
  1986.

\bibitem{BS-BTZ:97}
B.~Shapiro and B.~T. Zinn, ``High-frequency nonlinear vibrational control,''
  \emph{IEEE Transactions on Automatic Control}, vol.~42, no.~1, pp. 83--90,
  1997.

\bibitem{CX-TY-MI:2018}
X.~Cheng, Y.~Tan, and I.~Mareels, ``On robustness analysis of linear
  vibrational control systems,'' \emph{Automatica}, vol.~87, pp. 202--209,
  2018.

\bibitem{YQ-DSB-FP:22a}
Y.~Qin, D.~S. Bassett, and F.~Pasqualetti, ``Vibrational control of cluster
  synchronization: Connections with deep brain stimulation,'' in \emph{{IEEE}
  Conf.\ on Decision and Control}, Canc\'un, Mexico, Dec. 2022, to appear.

\bibitem{AC-MM:15}
A.~Chapman and M.~Mesbahi, ``State controllability, output controllability and
  stabilizability of networks: A symmetry perspective,'' in \emph{{IEEE} Conf.\
  on Decision and Control}.\hskip 1em plus 0.5em minus 0.4em\relax Osaka,
  Japan: IEEE, 2015, pp. 4776--4781.

\bibitem{TJ-THL-18}
J.~Trumpf and H.~L. Trentelman, ``Controllability and stabilizability of
  networks of linear systems,'' \emph{IEEE Transactions on Automatic Control},
  vol.~64, no.~8, pp. 3391--3398, 2018.

\bibitem{SP-SK-PAA:15}
S.~Pequito, S.~Kar, and P.~A. Aguiar, ``A framework for structural input/output
  and control configuration selection in large-scale systems,'' \emph{IEEE
  Transactions on Automatic Control}, vol.~61, no.~2, pp. 303--318, 2015.

\bibitem{LJ-CX-PS-PGJ-PVM:19}
J.~Li, X.~Chen, S.~Pequito, G.~J. Pappas, and V.~M. Preciado, ``Resilient
  structural stabilizability of undirected networks,'' in \emph{American
  Control Conference}, 2019, pp. 5173--5178.

\bibitem{DPC-TP:15}
C.~De~Persis and P.~Tesi, ``Input-to-state stabilizing control under
  denial-of-service,'' \emph{IEEE Transactions on Automatic Control}, vol.~60,
  no.~11, pp. 2930--2944, 2015.

\bibitem{DA-SF-DMD:16}
A.~D’Innocenzo, F.~Smarra, and M.~D. Di~Benedetto, ``Resilient stabilization
  of multi-hop control networks subject to malicious attacks,''
  \emph{Automatica}, vol.~71, pp. 1--9, 2016.

\bibitem{FP-SZ-FB:13q}
F.~Pasqualetti, S.~Zampieri, and F.~Bullo, ``Controllability metrics,
  limitations and algorithms for complex networks,'' \emph{IEEE Transactions on
  Control of Network Systems}, vol.~1, no.~1, pp. 40--52, 2014.

\bibitem{SSM-MH-MM:17}
S.~S. Mousavi, M.~Haeri, and M.~Mesbahi, ``On the structural and strong
  structural controllability of undirected networks,'' \emph{IEEE Transactions
  on Automatic Control}, vol.~63, no.~7, pp. 2234--2241, 2017.

\bibitem{JJ-HJW-HLT-KMC:19}
J.~Jia, H.~J. Van~Waarde, H.~L. Trentelman, and M.~K. Camlibel, ``A unifying
  framework for strong structural controllability,'' \emph{IEEE Transactions on
  Automatic Control}, vol.~66, no.~1, pp. 391--398, 2020.

\bibitem{GB-KO-AS-KM:21}
B.~Guo, O.~Karaca, S.~Azhdari, M.~Kamgarpour, and G.~Ferrari-Trecate,
  ``Actuator placement for structural controllability beyond strong
  connectivity and towards robustness,'' in \emph{{IEEE} Conf.\ on Decision and
  Control}, 2021, pp. 5294--5299.

\bibitem{AW-SM-YY-KX:22}
W.~Abbas, M.~Shabbir, Y.~Yaz{\i}c{\i}o{\u{g}}lu, and X.~Koutsoukos, ``Leader
  selection for strong structural controllability in networks using zero
  forcing sets,'' in \emph{American Control Conference}, 2022, pp. 1444--1449.

\bibitem{REB-JB-SMM:86a}
R.~E. Bellman, J.~Bentsman, and S.~M. Meerkov, ``Vibrational control of
  nonlinear systems: Vibrational stabilization,'' \emph{IEEE Transactions on
  Automatic Control}, vol.~31, no.~8, pp. 710--716, 1986.

\bibitem{HKK:02-bis}
H.~K. Khalil, \emph{Nonlinear Systems}.\hskip 1em plus 0.5em minus 0.4em\relax
  Prentice Hall, 2002.

\bibitem{QY-KY-BDOA-CM:2021}
Y.~Qin, Y.~Kawano, B.~D. Anderson, and M.~Cao, ``Partial exponential stability
  analysis of slow-fast systems via periodic averaging,'' \emph{IEEE
  Transactions on Automatic Control}, 2021.

\bibitem{DH-AJP:86}
D.~Hinrichsen and A.~J. Pritchard, ``Stability radii of linear systems,''
  \emph{Systems \& Control Letters}, vol.~7, no.~1, pp. 1--10, 1986.

\bibitem{LQ-BB-AR-EJD-PMY-JCD:95}
L.~Qiu, B.~Bernhardsson, A.~Rantzer, E.~J. Davison, P.~M. Young, and J.~C.
  Doyle, ``A formula for computation of the real stability radius,''
  \emph{Automatica}, vol.~31, no.~6, pp. 879--890, 1995.

\bibitem{JBJ-GG:00}
J.~Bang-Jensen and G.~Gutin, \emph{Digraphs: Theory, Algorithms and
  Applications}, ser. Monographs in Mathematics.\hskip 1em plus 0.5em minus
  0.4em\relax Springer, 2000.

\end{thebibliography}
	\newpage

\end{document}